\documentclass[12pt]{amsart}
\usepackage{amsfonts}
\usepackage{amssymb}
\theoremstyle{definition}

\theoremstyle{remark}

\newcommand{\ds}{\displaystyle}

\begin{document}

\centerline{\large\bf MATHEMATICS AND EDUCATION IN MATHEMATICS, 1988}
\centerline{\bf Proceedings of the Seventeenth Spring Conference of the Union of}
\centerline{\bf Bulgarian Mathematicians}
\centerline{\bf Sunny Beach, April 6-9, 1988}

\vspace{0.5in}
\centerline{\large ON THE EQUIVALENCE PROBLEM FOR MANIFOLDS}
\centerline{\large OF INDEFINITE METRIC }

\vspace{0.2in}
\centerline{\large  Ognian T. Kassabov }

\vspace{0.2in}
{\sl  Conditions, related to the Kulkarni's equivalence problem are considered for
indefinite Riemannian and Kaehlerian manifolds. Corresponding theorems are obtained
for the Ricci curvatures as well as for the holomorphic sectional curvatures of
indefinite Kaehler manifolds.}

\vspace{0.2in}

1. Let \,$M$\, be a Riemannian manifold of definite or indefinite metric \,$g$. A
2-plane \,$\alpha$\, in a tangent space is said to be nondegenerate, weakly
degenerate or strongly degenerate, if the rank of the restriction of \,$g$\, on\,
$\alpha$\, is 2, 1 or 0, respectively. A vector \,$\xi$\, on \,$M$\, is said to
be isotropic, if \,$g(\xi,\xi)=0$\, and \,$\xi\ne 0$. Of course, for degenerate
planes and isotropic vectors one speaks only when the metric is indefinite. The
curvature of a nondegenerate 2-plane \,$\alpha$\, with a basis \,$\{ x,y \}$\,
is defined by
$$
	K(\alpha) = \frac{R(x,y,y,x)}{\pi_1(x,y,y,x)} \ ,
$$
where \,$R$\, denotes the curvature tensor and
$$
	\pi_1(z,u,v,w)=g(z,w)g(u,v)-g(z,v)g(u,w) \ .
$$ 
Let \,$\overline M$\, be another Riemannian manifold of definite or indefinite metric. 
The correspond\-ing objects for \,$\overline M$\, will be denoted by a bar overhead. 
A diffeomorphism \,$f$\, of \,$M$\, onto \,$\overline M$\, is said to be sectional
curvature preserving [7], if
$$
	\overline K(f_*\alpha) = K(\alpha)    \leqno (1.1)
$$ 
for every nondegenerate 2-plane \,$\alpha$\, on \,$M$. In [7] R. S. Kulkarni investigates 
the converse of the so-called "theorema egregium" of Gauss and proves that any sectional 
curvature preserving diffeomorphism of Riemannian manifolds of nowhere constant sectional 
curva\-ture and dimension $\ge 4$ is trivial, i.e. it is an isometry. The related 
questions for the Ricci curvature and for the holomorphic sectional curvature of Kaehler
manifolds are considered in [8], [9]. Sectional curvature preserving diffeomorphisms of
indefinite Riemannian manifolds are studied in [10]. The condition, corresponding to
(1.1) for degenerate planes is
$$
	\lim_{\alpha\rightarrow\alpha_0} \frac{\overline K(f_*\alpha)}{K(\alpha)}=1 \ ,   \leqno (1.2)
$$
where the degenerate 2-plane \,$\alpha_0$\, is approximated by nondegenerate 2-planes, 
whose images are also nondegenerate. In [6] we
examine diffeomorphisms, satisfying (1.2) for weakly or for strongly degenerate 2-planes.
Then they appear the manifolds of quasi-constant curvature and the manifolds of recurrent 
curvature in the sense of Walker, whose definitions we recall:

An \,$n$-dimensional (indefinite) Riemannian manifold is said to be a \,$K_n^*$-manifold
[11], if either it is recurrent (i.e. \,$\nabla R=\alpha\otimes R$, where \,$\alpha \ne 0$) or
it is symmetric and there exists a differential form \,$\alpha\ne 0$, such that
$$
	\sum_{cycl\, x,y,z} \alpha(x)R(y,z,u,v)=0 \ .
$$
Walker [11] showed, that \,$\alpha$\, is defined by \,$\alpha(X)=g(\nabla v,X)$\,, where 
\,$v$\, is a smooth function (called recurrence-function) and \,$\nabla v$\, denotes
the gradient of \,$v$.

An \,$n$-dimensional (indefinite) Riemannian manifold is said to be of quasi-constant
curvature [1], [2], if it is conformally flat and there exist functions \,$H,\,N$\, 
and a unit vector \,$V$, such that the curvature tensor has the form
$$
	R=(N-H)\varphi(B)+H\pi_1 \ ,
$$
where \,$B(X,Y)=g(X,V)g(Y,V)$\, and \,$\varphi$\, is defined by
$$
	\varphi(Q)(x,y,z,u)=g(x,u)Q(y,z)-g(x,z)Q(y,u)+g(y,z)Q(x,u)-g(y,u)Q(x,z)
$$
for any symmetric tensor \,$Q$\, of type (0,2). Such a manifold we shall denote
by \,$M(H,N,V)$.

In this paper we consider diffeomorphisms of indefinite Riemannian manifolds,
satisfying conditions, analogous to (1.2) for the Ricci curvature and for the
holomorphic sectional curvature of indefinite Kaehler manifolds. We need
the following Lemma. 

\vspace{0.1in}
\underline{Lemma [6].} Let \,$f$\, be a diffeomorphism of indefinite Riemannian 
manifolds of dimension \,$n\ge 3$. Let in a point \,$p\in M$\, there exists an 
isotropic vector \,$\xi$, such that every isotropic vector which is sufficiently 
close to \,$\xi$\, is mapped by \,$f_*$\, in an isotropic vector in \,$f(p)$. 
Then \,$f$\, is a homothety in \,$p$.  

\vspace{0.2in}
\underline{2.} Let us recall that the Ricci curvature in the direction of a nonzero
nonisotropic vector \,$x$\, is defined by
$$
	K_S(x)=\frac{S(x,x)}{g(x,x)} \ ,
$$
where \,$S$\, is the Ricci tensor of \,$M$. As analogue of the Ricci curvature preserving 
diffeomorphisms [8], it is natural to consider diffeomorphisms, satisfying
$$
	\lim_{x\rightarrow\xi} \frac{\overline K_{\overline S}(f_*x)}{K_S(x)}=1 \ ,   \leqno (2.1)
$$
when the isotropic vector \,$\xi$\, is approximated by nonisotropic vectors, whose images 
are also nonisotropic. Then we have

\vspace{0.1in}
\underline{Theorem 1.} Let \,$M$\, and \,$\overline M$\, be indefinite Riemannian 
manifolds of dimension \,$n\ge 3$\, and let \,$f$\, be a diffeomorphism of \,$M$\,
onto \,$\overline M$\, satisfying (2.1) for every isotropic vector \,$\xi$\, on \,$M$.
If \,$M$\, is nowhere Einsteinian, then \,$f$\, is conformal. Let \,$f^*\bar g=e^{2\sigma}g$\,
and assume that \,$M$\, is conformally flat. Then:

a) if \,$\nabla\sigma$\, vanishes identically, then \,$f$\, is an isometry;

b) if \,$\nabla\sigma$\, is isotropic, then \,$M$\, is a conformal flat \,$K_n^*$\,-space
and \,$\sigma$\, is a function of the recurrence-function;

c) if \,$\parallel\nabla\sigma\parallel^2$\, doesn't vanish, then \,$M$\, is a manifold
\,$M(H,N,\nabla\sigma/\left\Arrowvert\nabla\sigma\right\Arrowvert)$\, of quasi-constant curvature,
\,$\nabla H$\, and \,$\nabla N$\, being proportional to \,$\nabla\sigma$.

\underline{ Proof.} \ Let \,$M$\, be non-Einsteinian in \,$p$, i.e. the Ricci tensor \,$S_p$\, is
not proportional to \,$g_p$. Then there exists an isotropic vector \,$\xi$\, in \,$p$, such that
\,$S(\xi,\xi)\ne 0$\, [4]. By continuity \,$S(\xi',\xi')\ne 0$\, for every isotropic vector
\,$\xi'$\,, sufficiently close to \,$\xi$. Then from (2.1) it follows that every such vector 
is mapped by \,$f_*$\, onto an isotropic one. By the Lemma \,$f$\, is a homothety in \,$p$.
Since the set of points in which \,$M$\, is not Einsteinian is dense, then \,$f$\, is
conformal.

Let \,$f^*\bar g=e^{2\sigma}g$. For the sake of simplicity of the denotations we identify
\,$M$\, with \,$\overline M$\, via \,$f$\, and omit \,$f_*$. Then (2.1) implies
\,$\overline S(\xi,\xi)=e^{2\sigma}S(\xi,\xi)$. Hence it follows (see [4])
$$
	\overline S=e^{2\sigma}\left\{ S+\frac{\bar\tau-\tau}ng \right\} \ ,  \leqno (2.2)
$$  
where \,$\tau$\, denotes the scalar curvature of \,$M$. Since \,$(M,\bar g)$\, and \,$(M,g)$\,
are conformally flat, their Weil conformal curvature tensors vanish [5], i.e.
$$
	\begin{array}{l}
		\ds \overline R=\frac1{n-2}\varphi(\overline S) - \frac{\bar\tau}{(n-1)(n-2)}\bar\pi_1 \ ,  \\
		\ds R=\frac1{n-2}\varphi(S) - \frac{\tau}{(n-1)(n-2)}\pi_1 \ .
	\end{array}  \leqno (2.3)
$$
From (2.2) and (2.3) it follows 
$$
	\overline R=e^{4\sigma} \left\{ R+\frac{\bar\tau - \tau}{n(n-1)}\pi_1 \right\} \ .
$$
Hence (1.2) is satisfied and the rest of the theorem follows from Theorem 2 in [6].

\vspace{0.2in}
\underline{3.} \ Let \,$M$\, be an indefinite Kaehler manifold with metric \,$g$\, and 
almost complex structure \,$J$. A 2-plane \,$\alpha$\, is said to be holomorphic, if 
\,$\alpha = J\alpha$. Note that a degenerate holomorphic 2-plane is necessarily
strongly degenerate. The Bochner curvature tensor \,$B$\, for \,$M$\, is defined by
$$
	B=R-\frac1{2(n+2)}(\varphi+\psi)(S)+\frac{\tau}{4(n+1)(n+2)}(\pi_1+\pi_2) \ ,
$$
where \,$2n$\, is the dimension of \,$M$,
$$
	\begin{array}{r}
		\psi(Q)(x,y,z,u)=g(x,Ju)Q(y,Jz)-g(x,Jz)Q(y,Ju)-2g(x,Jy)Q(z,Ju)  \\
		                +g(y,Jz)Q(x,Ju)-g(y,Ju)Q(x,Jz)-2g(z,Ju)Q(x,Jy)
	\end{array}  
$$
for a symmetric tensor \,$Q$\, of type (0,2) and \,$\pi_2=\frac12\psi(g)$. As usual we 
denote the curvature of a nondegenerate holomorphic 2-plane with an orthonormal basis
\,$\{x,Jx\}$\, by \,$H(x)$. For the holomorphic curvatures, i.e. the curvatures of
holomorphic planes, the analogue of (1.2) is
$$
	\lim_{x\rightarrow\xi} \frac{\overline H(f_*x)}{H(x)}=1 \ ,   \leqno (3.1)
$$
where the isotropic vector \,$\xi$\, is approximated by nonisotropic vectors, whose
images are also nonisotripic and with the natural requirement that the image of any
holomorphic 2-plane is also holomorphic (see also [9]). Then we have

\vspace{0.1in}
\underline{Theorem 2.} Let \,$M$\, and \,$\overline M$\, be indefinite Kaehler manifolds
of dimension \,$2n \ge 4$\, and let \,$f$\, be a diffeomorphism of \,$M$\, onto 
\,$\overline M$\,, satisfying (3.1) for every isotropic vector \,$\xi$\, on \,$M$. If
\,$M$\, is not of constant holomorphic sectional curvature, then \,$f$\, is a 
holomorphic or antiholomorphic isometry.

\underline{ Proof.} \ By Lemma 2 in [9] \,$f$\, is holomorphic or antiholomorphic. Let
\,$N$\, be the set of points, in which \,$M$\, is not of constant holomorphic sectional 
curvature and let \,$p\in N$. To show that \,$f$\, is a homothety in \,$p$, we shall 
consider two cases:

1) The Bochner curvature tensor of \,$M$\, vanishes in \,$p$. Then \,$M$\, cannot be Einsteinian
in \,$p$, because it is not of constant holomorphic sectional curvature in \,$p$. 
Hence there exists an isotropic vector \,$\xi$\, in \,$T_pM$\,, such that 
\,$S(\xi,\xi)\ne 0$. Then \,$S(\xi',\xi')\ne 0$\, for every isotropic vector \,$\xi'$,
sufficiently close to \,$\xi$. By (3.1) and \,$B=0$\, we have
$$
	1=\lim_{x\rightarrow\xi'} \frac{\overline H(f_*x)}{H(x)} =
		\lim_{x\rightarrow\xi'} \frac{\overline H(f_*x)g(x,x)}{\frac4{n+2}S(x,x)-
						\frac{\tau}{(n+1)(n+2)}g(x,x)}
$$
$$
	=\frac{n+2}{4S(\xi',\xi')}\lim_{x\rightarrow\xi'} 
			\frac{\overline R(f_*x,Jf_*x,Jf_*x,f_*x)g(x,x)}{\bar g^2(f_*x,f_*x)} \ .
$$
Hence \,$\bar g(f_*\xi',f_*\xi')=0$, i.e. \,$f_*\xi'$\, is isotropic. According to the
Lemma from section 1 \,$f$\, is a homothety in \,$p$. 

2) The Bochner curvature tensor doesn't vanish in \,$p$. Then there exists an isotropic 
vector \,$\xi$\, in \,$T_pM$, such that \,$R(\xi,J\xi,J\xi,\xi)\ne 0$ [3]. Consequently
\,$R(\xi',J\xi',J\xi',\xi')\ne 0$\, for every isotropic vector \,$\xi'$, sufficiently
close to \,$\xi$. Then
$$
	1=\lim_{x\rightarrow\xi'} \frac{\overline H(f_*x)}{H(x)} =
		\frac1{R(\xi',J\xi',J\xi',\xi')}\,\lim_{x\rightarrow\xi'} 
		\frac{\overline R(f_*x,Jf_*x,Jf_*x,f_*x)}{\bar g^2(f_*x,f_*x)}g^2(x,x) 
$$ 
which implies again \,$\bar g(f_*\xi',f_*\xi')=0$\, and hence \,$f$\, is a
homothety in \,$p$. 

So \,$f$\, is a homothety in \,$p$. Since \,$M$\, is not of constant holomorphic 
sectional curvature, the set \,$N$\, is dense in \,$M$. This implies that 
\,$f$\, is conformal, i.e. \,$f^*g=\lambda g$, where \,$\lambda$\, is a smooth
function. Denote by \,$\Phi$\, the fundamental form of \,$M$ : \ $\Phi(X,Y)=g(JX,Y)$.
Then \,$f^*\overline\Phi=\lambda\Phi$. Since \,$\Phi$\, and $\overline\Phi$\, are
closed, this yelds \,$\alpha\lambda \wedge\Phi=0$\, (see also [9]), which implies
that \,$\lambda$\, is a (nonzero) constant, thus proving the theorem. 

\vspace{0.1in}
\underline{Corollary.} \ If in Theorem 2 \,$M$\, has nonvanishing Bochner curvature 
tensor, then \,$f$\, is a holomorphic or antiholomorphic isometry.

\underline{Proof.} \,From \,$f^*\bar g=\lambda R$\, it follows 
\,$f^*\overline R=\lambda R$. Hence
$$
	(f^*\overline R)(\xi,J\xi,J\xi,\xi)=\lambda R(\xi,J\xi,J\xi,\xi) \ .   \leqno (3.2)
$$
On the other hand (3.1) implies
$$
	(f^*\overline R)(\xi,J\xi,J\xi,\xi)=\lambda^2R(\xi,J\xi,J\xi,\xi)     \leqno (3.3)
$$
for every isotropic vector \,$\xi$\, on \,$M$. Since \,$M$\, has nonvanishing
Bochner curvature tensor, there exists a point \,$p$\, in \,$M$\, and an
isotropic vector \,$\xi$\, in \,$T_pM$, such that $R(\xi,J\xi,J\xi,\xi)\ne 0$\,
[3]. Then (3.2) and (3.3) yield \,$\lambda=1$, proving the assertion.

\vspace{0.1in}
Note that if (2.1) is satisfied for every isotropic vector for indefinite
Kaehler manifold and \,$M$\, is not Einsteinian, as in Theorem 1 we obtain 
\,$f^*\bar g=\lambda g$\, for $\lambda \in \mathfrak FM$. Assume that \,$f$\, is holomorphic or 
antiholomorphic. Then as in Theorem 2 \,$\lambda$\, is a constant and similarly
to the case in the Corollary \,$\lambda =1$. Thus we have

\vspace{0.1in}
\underline{Theorem 3.} Let \,$M$\, and \,$\overline M$\, be indefinite Kaehler manifolds
of dimension \,$2n \ge 4$\, and let \,$f$\, be a holomorphic or antiholomorphic 
diffeomorphism of \,$M$\, onto 
\,$\overline M$\,, satisfying (2.1) for every isotropic vector \,$\xi$\, on \,$M$. If
\,$M$\, is not Einsteinian, then \,$f$\, is an isometry.

The same assertion holds for diffeomorphisms of definite or indefinite Kaehler
manifolds if (2.1) is changed by
$$
	\overline K_{\overline S}(f_*x)=K_S(x)
$$
for every unit vector \,$x$\, on \,$M$.

\vspace{0.5in}
\centerline{\large REFERENCES}

\vspace{0.1in}
\noindent
1. T. Adati, Y. Wong. Manifolds of quasi-constant curvature I: A manifold
of quasi-

constant curvature and an $s$-manifold. TRU. Math., 21(1985), 95-103.

\noindent
2. V. Boju, M. Popescu: Espaces \`a courbure quasi-constante. J. Differ. Geom., 13(1978), 

375-383.

\noindent
3. A. Borisov, G. Ganchev, O. Kassabov. Curvature properties and isotropic planes
of 

Riemannian and almost Hermitian manifolds of indefinite metrics.
Ann. Univ. Sof., 

Fac. Math.  M\'ec., 78(1984), 121-131.

\noindent
4. M. Dajczer, K. Nomizu. On boundedness of Ricci curvature of an indefinite metric.

Bol. Soc. Brasil. Mat., 11(1980), 25-30. 

\noindent
5. L. P. Eisenhart. Riemannian geometry. Princeton. University Press, 1949.

\noindent
6. O. Kassabov. Diffeomorphisms of pseudo-Riemannian manifolds and the values of the

curvature tensor on degenerate planes. Serdica, 15(1989), 78-86.
 
\noindent
7. R. S. Kulkarni. Curvature and metric. Ann. of Math., 91(1970), 311-331.

\noindent
8. R. S. Kulkarni. Curvature structures and conformal transformations. J. Differ. Geom.,

4(1970), 425-451.

\noindent
9. R. S. Kulkarni. Equivalence of Kaehler manifolds and other equivalence problems. J.

Differ. Geom., 9(1974), 401-408. 

\noindent
10. B. Ruh. Krummungstreue Diffeomorphismen Riemannscher und pseudo-Riemannscher

Mannigfaltigkieten. Math. Z., 189(1985), 371-391.

\noindent
11. A. G. Walker. On Ruse's spaces of recurrent curvature. Proc. Lond. Math. Soc., 
II 

S\'er., 52(1951), 36-64.

\end{document}